\documentclass{article}
\usepackage{amsmath}
\usepackage{amsfonts}
\usepackage{amsthm}
\usepackage{amssymb}

\usepackage{centernot}
\usepackage{url}

\def\mod{\bmod}

\def\ratQ{ \mathbb{Q} }  
 
\def\natN{  \mathbb{N} }

\numberwithin{equation}{section}

\newtheorem{Theorem}{Theorem}[section]
\newtheorem*{theorem*}{Theorem 1.2}
\newtheorem{Corollary}[Theorem]{Corollary}

\newtheorem{Lemma}[Theorem]{Lemma}

\newtheorem{Conjecture}[Theorem]{Conjecture}

\begin{document}

On a conjecture concerning the number of solutions to $a^x+b^y=c^z$  

13 July 2022      

Maohua Le, Reese Scott, Robert Styer

\

\begin{abstract}
Let $a$, $b$, $c$ be fixed coprime positive integers with $\min\{ a,b,c \} >1$.  Let $N(a,b,c)$ denote the number of positive integer solutions $(x,y,z)$ of the equation $a^x + b^y = c^z$.  We show that if $(a,b,c)$ is a triple of distinct primes for which $N(a,b,c)>1$ and $(a,b,c)$ is not one of the six known such triples then $c>10^{18}$, and there are exactly two solutions $(x_1, y_1, z_1)$, $(x_2, y_2, z_2)$ with $2 \mid x_1$, $2 \mid y_1$, $z_1=1$, $2 \nmid y_2$, $z_2>1$, and, taking $a<b$, we must have $a=2$, $b \equiv 1 \bmod 12$, $c \equiv 5 \mod 12$, with $(a,b,c)$ satisfying further strong restrictions.  These results support a conjecture put forward in \cite{ScSt6}.  
\end{abstract}

2020 Mathematics Subject Classification 11D61


\section{Introduction}  %1

Let $\natN$ be the set of all positive integers. Let $a$, $b$, $c$ be fixed coprime positive integers with $\min\{ a,b,c\}>1$.  The equation
$$a^x + b^y = c^z, x,y,z \in \natN \eqno{(1.1)}$$
has been studied deeply with abundant results (see \cite{LeScSt}).  In 1933, K. Mahler \cite{M} used his $p$-adic analogue of the Diophantine approximation method of Thue-Siegel to prove that (1.1) has only finitely many solutions $(x,y,z)$.  However, his method is ineffective.  Let $N(a,b,c)$ denote the number of solutions $(x,y,z)$ of (1.1).  An effective upper bound for $N(a,b,c)$ was first given by A. O. Gel'fond \cite{G}.  A straightforward application of an upper bound on the number of solutions of binary $S$-unit equations due to F. Beukers and H. P. Schlickewei \cite{BSch} gives $N(a,b,c) \le 2^{36}$.  The following more accurate upper bounds for $N(a,b,c)$ have been obtained in recent years:

(i)  (R. Scott and R. Styer \cite{ScSt6})  If $ 2 \nmid c$ then $N(a,b,c) \le 2$.

(ii)  (Y. Z. Hu and M. H. Le \cite{HL1})  If $\max\{ a,b,c\} > 5 \cdot 10^{27}$, then $N(a,b,c) \le 3$.  

(iii)   (Y. Z. Hu and M. H. Le \cite{HL2})  If $2 \mid c$ and $\max\{ a,b,c\} > 10^{62}$, then $N(a,b,c) \le 2$.

(iv)  (T. Miyazaki and I. Pink \cite{MP})  If $2 \mid c$, $a<b$, and $\max\{a,b,c \} \le 10^{62}$, then $N(a,b,c) \le 2$ except for $N(3,5,2) = 3$.  

Nevertheless, the problem of establishing $N(a,b,c) \le 1$ with a finite number of exceptions remains open.  This open question is addressed by the following conjecture in \cite{ScSt6}.  Assuming without loss of generality that $a$, $b$, and $c$ are not perfect powers, the conjecture may be formulated as follows: 

\begin{Conjecture}  % 1.1  
For $a<b$, we have $N(a,b,c) \le 1$, except for 

(i)  $N(2, 2^r-1, 2^r+1)=2$, $(x,y,z)=(1,1,1)$ and $(r+2, 2, 2)$, where $r$ is a positive integer with $r \ge 2$.  

(ii)  $N(2, 3, 11)=2$, $(x,y,z)=(1,2,1)$ and $(3,1,1)$.

(iii)  $N(2, 3, 35)=2$, $(x,y,z)=(3,3,1)$ and $(5,1,1)$.

(iv)  $N(2, 3, 259)=2$, $(x,y,z)=(4,5,1)$ and $(8,1,1)$.

(v)  $N(2, 5, 3)=2$, $(x,y,z)=(1,2,3)$ and $(2,1,2)$.

(vi)  $N(2, 5, 133)=2$, $(x,y,z)=(3,3,1)$ and $(7,1,1)$.

(vii)  $N(2, 7, 3)=2$, $(x,y,z)=(1,1,2)$ and $(5,2,4)$.

(viii)  $N(2, 89, 91)=2$, $(x,y,z)=(1,1,1)$ and $(13,1,2)$.

(ix)  $N(2, 91, 8283)=2$, $(x,y,z)=(1,2,1)$ and $(13,1,1)$.

(x)  $N(3,5,2)=3$, $(x,y,z)=(1,1,3)$, $(1,3,7)$, and $(3,1,5)$.

(xi)  $N(3,10,13)=2$, $(x,y,z)=(1,1,1)$ and $(7,1,3)$.

(xii)  $N(3,13,2)=2$, $(x,y,z)=(1,1,4)$ and $(5,1,8)$.

(xiii)  $N(3, 13, 2200)=2$, $(x,y,z)=(1,3,1)$ and $(7,1,1)$.

\end{Conjecture}

Later in this paper, in referring to the solutions in (i) through (xiii) above, it will be helpful to have established the following result:

\begin{Lemma} %1.2
The values of $N(a,b,c)$ in Conjecture 1.1 are exact: there are no further solutions in any of the thirteen cases.  
\end{Lemma}  

\begin{proof} 
For cases (i) through (xii) this follows from Theorem 1 of \cite{ScSt6} and Theorem 6 of \cite{Sc}.  For (xiii), consideration modulo 16 shows that $z>1$ requires $4 \mid x-y$, contradicting consideration modulo 5.  So $z=1$, and there are only two solutions.  
\end{proof}

Although the results on (1.1) known so far support Conjecture 1.1, it is generally far from being resolved.  This difficult problem is made more approachable by taking $c$ prime.  More than thirty years ago, the first author \cite{Le1} discussed the upper bound for $N(a,b,c)$ when $a$, $b$, $c$ are distinct primes.  Many authors have used this approach to the problem.  Later, \cite{ScSt6} removed the difficulty caused by taking $c$ composite when $c$ is odd, and Hu and Le \cite{HL2} and Miyazaki and Pink \cite{MP} handled even composite $c$; these later results established $N(a,b,c) \le 2$ with the single exceptional case $(a,b,c) = (3,5,2)$.  

Establishing $N(a,b,c) \le 1$ involves many exceptional cases and is in general much more difficult, suggesting perhaps that the old approach of considering only prime bases may be a practical way to begin considering this problem.  

Our first result requires $c$ prime with congruence restrictions on $a$ and $b$ not necessarily prime: 

\begin{Theorem} % 1.3  
If $a \equiv 2 \bmod 3$, $b \not\equiv 0 \bmod 3$, and $c>3$ is an odd prime, then $N(a,b,c)  \le 1$, except for the following possibility: 
(1.1) has exactly two solutions $(x_1, y_1, z_1)$ and $(x_2, y_2, z_2)$, and these solutions satisfy $2 \mid x_1$, $2 \mid y_1$, $z_1=1$, and $2 \nmid y_2$, with $c \equiv 5 \bmod 12$.  
\end{Theorem} 

Theorem 1.3 is used to establish the following result: 

\begin{Theorem} %1.4 
If $a$, $b$, $c$ are distinct primes with $a<b$ and $N(a,b,c) > 1$, then $c \equiv 5 \bmod 12$, and there must be exactly two solutions $(x_1,y_1,z_1)$ and $(x_2,y_2,z_2)$, with $2 \mid x_1$, $2 \mid y_1$, $z_1=1$, and $2 \nmid y_2$, except for $(a,b,c) = (2,3,5)$, $(2,3,11)$, $(2,5,3)$, $(2,7,3)$, $(3,5,2)$, and $(3,13,2)$. 
\end{Theorem}  

Further restrictions are given by the following: 

\begin{Theorem}  %1.5
If $a$, $b$, $c$ are distinct primes with $a<b$ and $N(a,b,c) > 1$, then, if $(a,b,c) \ne (2,3,5)$, $(2,3,11)$, $(2,5,3)$, $(2,7,3)$, $(3,5,2)$, or $(3,13,2)$, we must have $a=2$ with $b \equiv 1 \bmod 12$ and $c \equiv 5 \bmod 12$; further, if $b \equiv 1 \bmod 24$, then $c \equiv 17 \bmod 24$.      
\end{Theorem}

\begin{Corollary}  % 1.6
If $a$, $b$, $c$ are distinct primes with $a<b$ and $N(a,b,c)>1$, then, if $(a,b,c) \ne (2,3,5)$, $(2,3,11)$, $(2,5,3)$, $(2,7,3)$, $(3,5,2)$, or $(3,13,2)$, we must have $N(a,b,c) = 2$ and, letting the two solutions be $(x_1,y_1,z_1)$ and $(x_2,y_2,z_2)$ as in Theorem 1.4, we must have
$$ z_2 > 1.  \eqno{(1.2)}$$
\end{Corollary}

In a later section we will use the following version of Conjecture 1.1 in which $a$, $b$, and $c$ are restricted to prime values: 

\begin{Conjecture}  % 1.7  
For $a$, $b$, and $c$ distinct primes with $a<b$, we have $N(a,b,c) \le 1$, except for 

(i)  $N(2,3,5) = 2$, $(x,y,z) = (1,1,1)$ and $(4,2,2)$.  

(ii)  $N(2, 3, 11)=2$, $(x,y,z)=(1,2,1)$ and $(3,1,1)$.

(iii)  $N(2, 5, 3)=2$, $(x,y,z)=(1,2,3)$ and $(2,1,2)$.

(iv)  $N(2, 7, 3)=2$, $(x,y,z)=(1,1,2)$ and $(5,2,4)$.

(v)  $N(3,5,2)=3$, $(x,y,z)=(1,1,3)$, $(1,3,7)$, and $(3,1,5)$.

(vi)  $N(3,13,2)=2$, $(x,y,z)=(1,1,4)$ and $(5,1,8)$.

\end{Conjecture}

In the final section of this paper, we will explain a method by which we have shown

\begin{Theorem}   % 1.8
Any counterexample to Conjecture 1.7 must have
$$ b > 10^9, c>10^{18}.$$
\end{Theorem}

This result is in marked contrast to results which can be obtained without assuming $a$, $b$, and $c$ are prime: in that more general case, the lower bound on $c$ when $(a,b,c)$ gives a counterexample to Conjecture 1.1 is still quite low (the latest such results are given by Miyazaki and Pink \cite{MP2} in considering Conjecture 1.1 for the special case in which both $a$ and $b$ are congruent to $\pm 1$ modulo $c$).  Theorem 1.8 improves Lemma 2.9 in Section 2 which follows.

\section{Preliminaries}  %2

We now divide all solutions $(x,y,z)$ of (1.1) into four classes according to the parities of $x$ and $y$: $2\mid x$ and $2 \mid y$, $2 \nmid x$ and $2 \mid y$, $2 \mid x$ and $2 \nmid y$, or $2 \nmid x$ and $2 \nmid y$.  We will call these classes the {\it parity classes} of (1.1).

\begin{Lemma}  % 2.1
If $c$ is an odd prime, then, for a given parity class, there is at most one solution $(x,y,z)$ to (1.1), except for $(a,b,c)=(3,10,13)$ or $(10,3,13)$.  
\end{Lemma}

\begin{proof} 
Since $c$ is an odd prime, using the notation of \cite{ScSt6} we see that for any given parity class of (1.1) there is only one ideal factorization in $C_D$.  Therefore, by Lemma 2 of \cite{ScSt6}, we obtain the lemma immediately. 
\end{proof}  

\begin{Lemma}[Lemma 2 of \cite{ScSt1}] % 2.2
The equation 
$$3^x+2^y=n^z \eqno{(2.1)}$$
has no solutions in positive integers $(x,y,z,n)$ with $z>1$ except for $3^2+2^4 = 5^2$.  
\end{Lemma} 

\begin{Lemma} %2.3
Let $X$, $n$ be positive integers with $2 \nmid X$ and $n>1$.  Then we have $\mid X^2 - 2^n \mid > 2^{0.26 n}$, except for $X^2 - 2^n = 1$ or $-7$.  
\end{Lemma}

\begin{proof} 
This lemma is a special case of Corollary 1.7 of \cite{BB} with $y=2$.  
\end{proof}

\begin{Lemma}[\cite{Co}, \cite{Le}] %2.4 
The equation
$$ X^2 + 2^m = Y^n, X,Y,m,n \in \natN, \gcd(X,Y)=1, n>2, \eqno{(2.2)}$$
has only the solutions $(X,Y,m,n) = (5,3,1,3)$, $(11,5,2,3)$, and $(7,3,5,4)$.  
\end{Lemma}

\begin{Lemma}[Theorem 8.4 of \cite{BS}]  % 2.5
The equation 
$$ X^2 - 2^m = Y^n, X,Y,m,n \in \natN, \gcd(X,Y)=1, m>1, n>2, \eqno{(2.3)}$$
has only the solution $(X,Y,m,n) = (71, 17, 7, 3)$.  
\end{Lemma}

\begin{Lemma}[Theorem 6 of \cite{ScSt1}]  %2.6
Let $A$, $B$ be distinct odd positive primes.  For a given positive integer $k$, the equation
$$ A^m - B^n = 2^k, m, n \in \natN, \eqno{(2.4)} $$
has at most one solution in positive integers $(m,n)$.  
\end{Lemma}

\begin{Lemma}[Theorem 6 of \cite{Sc}]  %2.7 
If $a<b$ and $c=2$, then $N(a,b,2) \le 1$, except for $N(3,5,2)=3$ and $N(3,13,2)=2$.  
\end{Lemma}

\begin{Lemma}[\cite{ScSt6}]  %2.8
If $2 \nmid c$, then $N(a,b,c) \le 2$. 
\end{Lemma}

\begin{Lemma}[\cite{Cao}]  %2.9
If $a$, $b$, $c$ are distinct primes with $\max\{ a,b,c\} < 100$, then Conjecture 1.7 is true.  
\end{Lemma}

\begin{Lemma}[Lemma 4.2 of \cite{ScSt7}] %2.10
The equation   
$$ (1+ \sqrt{-D} )^r = m \pm \sqrt{-D} \eqno{(2.5)} $$
has no solutions with integer $r>1$ when $D$ is a positive integer congruent to 2 mod 4 and $m$ is any integer, except for $D=2$, $r=3$.    
  
Further, when $D \equiv 0 \bmod 4$ is a positive integer such that $1+D$ is prime or a prime power, (2.5) has no solutions with integer $r>1$ except for $D=4$, $r=3$.    
\end{Lemma}

\begin{Lemma}  %2.11
If the equation 
$$ 3^m - 2^n = 3^x - 2^y = d \eqno{(2.6)}$$
has a solution in positive integers $(m,n,x,y)$ with $m \ne x$, then $d=1$, $-5$, or $-13$.  
\end{Lemma}

\begin{proof} 
This follows easily from a conjecture of Pillai \cite{P} first proven by Stroeker and Tijdeman \cite{StTi} using lower bonds on linear forms in logarithms.  It is also an immediate consequence of the elementary Corollary to Theorem 4 of \cite{Sc}.   
\end{proof}

\begin{Lemma} %  2.12
If $N(a,b,c)>1$ when $a \equiv 2 \bmod 3$, $b \not\equiv 0 \mod 3$, and $c$ is an odd prime,  then (1.1) has exactly two solutions $(x_1, y_1, z_1)$ and $(x_2,y_2,z_2)$, where $2 \nmid y_1 - y_2$. 
\end{Lemma} 

\begin{proof}
If $a \equiv 2 \bmod 3$, $b \not\equiv  0 \bmod 3$, and $c$ is an odd prime, then the parity of $y$ determines the parity class, so that, since Lemma 2.1 shows that there is at most one solution per parity class, we see that $N(a,b,c)>1$ implies that there are exactly two solutions $(x_1, y_1, z_1)$ and $(x_2,y_2,z_2)$ and $ 2 \nmid y_1 - y_2$. 
\end{proof} 

\begin{Lemma}[Theorem 1.1 of \cite{Be2}]  % 2.13
Let $c$ and $b$ be positive integers.  Then there exists at most one pair of positive integers $(z,y)$ for which 
$$ 0 < |c^{z} - b^y| < \frac{1}{4} \max\{ c^{z/2}, b^{y/2} \}.$$
\end{Lemma}

\section{Proofs of Theorems 1.3 and 1.4}  %3 

We prove Theorem 1.3:

\begin{proof}  
Assume $N(a,b,c)>1$ with $a \equiv 2 \bmod 3$, $b \not\equiv 0 \mod 3$, and $c>3$ an odd prime.  By Lemma 2.12 there must be exactly two solutions $(x_1, y_1, z_1)$ and $(x_2,y_2,z_2)$, with $2 \nmid y_1 - y_2$.  Take $y_1$ even.  Then, since $c>3$, consideration modulo 3 shows that $2 \mid x_1$, $2 \mid y_1$, and $2 \nmid z_1$, which requires $c \equiv 2 \bmod 3$ and $c \equiv 1 \bmod 4$.  Thus, $c \equiv 5 \bmod 12$.  To see that $z_1=1$, note that, since $c>3$ and $2 \nmid z_1$, it suffices to use Lemma 2.4, noting $(a,b,c) \ne (2,5,3)$ or $(2,7,3)$ and handling $(2,11,5)$ using Lemma 2.1 with consideration modulo 3 and modulo 5.   
\end{proof} 

We now prove Theorem 1.4:

\begin{proof}
Assume that $a$, $b$, and $c$ are distinct primes with $a<b$.  Then, if $N(a,b,c)>1$, we can take $a=2$ ($c \ne 2$ by Lemma 2.7, except for $(a,b,c)=(3,5,2)$ and $(3,13,2)$, given as exceptions in the statement of the theorem).  

And we can also take $b \ne 3$: using Lemma 2.2, we find that if the equation $2^x+3^y=c^z$ has two solutions $(x_1,y_1,z_1)$ and $(x_2,y_2,z_2)$, we either must have $c=5$ with $\{z_1,z_2 \}=\{ 1,2 \}$, or we must have $z_1=z_2=1$, in which case we have $3^{y_1}-2^{x_2} = 3^{y_2}-2^{x_1} = d$, where $d=1$, $-5$, or $-13$ by Lemma 2.11; recalling Lemma 1.2 and using the solutions given in (x) and (iii) in Conjecture 1.1, we see that  $d = -5$ corresponds to $c=35$, while using (xii) with (iv) shows that $d=-13$ corresponds to $c=259$; since neither of these values of $c$ is prime we must have $d=1$, and it is a familiar elementary result that we must have $(x_1, y_1, x_2, y_2) = (1,2,3,1)$ giving $c=11$. 

Thus we find that the only cases with $a=2$, $b=3$, and $N(a,b,c)>1$ are $(a,b,c) = (2,3,5)$ and $(2,3,11)$, given as exceptions in the statement of the theorem.     

Now assume $a=2$ and $c=3$, and assume (1.1) has two solutions $(x_1,y_1,z_1)$ and $(x_2,y_2,z_2)$.  By Lemma 2.12, we can assume $2 \mid y_1$.  If $z_1$ is even, then $3^{z_1/2} - b^{y_1/2} = 2$ and $3^{z_1/2}+b^{y_1/2} = 2^{x_1-1}$, so that $ 3^{z_1/2} = 2^{x_1 -2} +1$; it is a familiar elementary result that we must have either $(x_1, z_1) = (3,2)$ (giving $b^{y_1/2} = 1$) or $(x_1, z_1) = (5,4)$ (giving $b^{y_1/2} = 7$).  Since $b>1$ we find that 
$(a,b,c) = (2,7,3)$ is the only possibility when $z_1$ is even.  

And if $2 \nmid z_1$ when $a=2$ and $c=3$, then, since $2 \mid y_1$, we have $x_1 = 1$, and we can factor in $\ratQ(\sqrt{-2})$: 
$$ \pm b^{y_1/2} \pm \sqrt{-2} = (1+\sqrt{-2})^{z_1}  \eqno{(3.1)} $$
where the \lq $\pm$' are independent.  Since clearly $z_1>1$, we can use Lemma 2.10 to see that $z_1 =3$ and $(a,b,c)=(2,5,3)$.  

Thus, recalling Lemma 1.2, we find that the only cases with $a=2$, $c=3$, and $N(a,b,c)>1$ are $(a,b,c)=(2,7,3)$ and $(2,5,3)$, given as exceptions in the statement of the theorem. 

So, excluding the six exceptions given in the statement of the theorem, we can assume $a=2$, $b \ne 3$, and $c \ne 3$ when $N(a,b,c)>1$ for a triple of primes $(a,b,c)$.  Now Theorem 1.4 follows immediately from Theorem 1.3.  
\end{proof}

\section{Proofs of Theorem 1.5 and Corollary 1.6}     %4

We now prove Theorem 1.5:

\begin{proof}
Assume $a$, $b$, and $c$ are distinct primes with $a<b$, $(a,b,c) \ne (2,3,5)$, $(2,3,11)$, $(2,5,3)$, $(2,7,3)$, $(3,5,2)$, or $(3,13,2)$, and $N(a,b,c) > 1$.  
From Theorem 1.4 we have 
$$ c \equiv 5 \bmod 12 \eqno{(4.1)}$$
and 
$$ 2^{x_1} + b^{y_1} = c, 2 \mid x_1, 2 \mid y_1.  \eqno{(4.2)}$$

Assume first 
$$b \equiv 1 \bmod 3. \eqno{(4.3)}$$
We have 
$$2^{x_2} + b^{y_2} = c^{z_2}, 2 \mid x_2, 2 \nmid y_2, 2 \nmid z_2.  \eqno{(4.4)}$$
Since $2 \mid x_2$ and $2 \nmid y_2$, by (4.1) and (4.4), we get $ b \equiv b^{y_2} \equiv c^{z_2} - 2^{x_2} \equiv 1-0 \equiv 1 \bmod 4$.  Hence, by (4.3), we have 
$$ b \equiv 1 \bmod 12 \eqno{(4.5)}$$
If $c \equiv 5 \bmod 8$ and $b \equiv 1 \bmod 8$, then from (4.2) and (4.4) we get $x_1=2$ and $x_2=2$, respectively.  This implies that (2.4) has two solutions $(m,n)= (1, y_1)$ and $(z_2, y_2)$  for $(A,B) = (c,b)$ and $k=2$.  By Lemma 2.6 this is impossible.  Hence, by (4.1) and (4.5), we have $(b,c) \not\equiv (1,5) \bmod 24$ and 
$$ (b,c) \equiv (1,17), (13,5), {\rm \ or \ } (13,17) \bmod 24. \eqno{(4.6)}$$
Now assume $ b \equiv 2 \bmod 3$.  Then we have (4.1), (4.2), and 
$$ 2^{x_2} + b^{y_2} = c^{z_2}, 2 \nmid x_2, 2 \nmid y_2, \eqno{(4.7)}$$
$$ b \equiv 2 \bmod 3, \eqno{(4.8)}$$
and
$$ c^{z_2} \equiv 1 \bmod 3.  \eqno{(4.9)}$$
Further, by (4.1) and (4.9), we get
$$ 2 \mid z_2.  \eqno{(4.10)}$$

When $x_2>1$ and $y_2>1$, we see from (4.7) and (4.10) that (2.3) has a solution $(X,Y,m,n)=(c^{z_2/2}, b, x_2, y_2)$.  Hence, by Lemma 2.5, we get $(b,c) = (17,71)$, for which a solution $(x_1, y_1, z_1)$ is impossible since $z_1=1$.  So we have either $x_2=1$ or $y_2=1$.  

If $y_2=1$, then from (4.7) and (4.10) we get
$$ 2^{x_2} + b = c^{z_2}, 2 \nmid x_2, 2 \mid z_2.  \eqno{(4.11)}$$
Note that $x_2>1$ since $c^{z_2}>c$.  

Apply Lemma 2.3 to (4.11) to obtain
$$ b = \left( c^{z_2/2} \right)^2 - 2^{x_2} > 2^{0.26 x_2}. \eqno{(4.12)}$$
Further, by (4.2), (4.11), and (4.12), we get
$$ 2^{x_2} + b = c^{z_2} \ge c^2 = (2^{x_1} + b^{y_1})^2 > b^{2y_1} \ge b^4 > 2^{1.04 x_2}, \eqno{(4.13)}$$
whence we obtain
$$ b > 2^{x_2} (2^{0.04 x_2} -1).  \eqno{(4.14)}$$
On the other hand, by (4.2), we have $b < \sqrt{c}$.  Hence, since by Lemma 2.9 we can assume $\max\{ b,c\} > 100$, by (4.11) we have $2^{x_2} = c^{z_2} - b > c^2 - \sqrt{c} \ge 101^2 - \sqrt{101} > 10190$, whence we get $x_2 \ge 15$.  So we have $2^{0.04 x_2} - 1 \ge 2^{0.6}-1 > 1/2$.  By (4.14), we get 
$$ b > 2^{x_2-1}.  \eqno{(4.15)}$$
Therefore, by (4.11), (4.13), and (4.15), we have $3b > 2^{x_2} + b = c^{z_2} \ge c^2 > b^4$, a contradiction.   So we obtain $y_2>1$.  This implies 
$$ x_2=1. \eqno{(4.16)}$$  
By (4.16), (4.1), (4.7), and (4.8) we have 
$$ b = 2^r h - 1, c = 2^s k + 1, r >1, s>1, 2 \nmid h, 2 \nmid k, 3 \mid h. \eqno{(4.17)}$$
Recalling (4.16), (4.7), and (4.10), we have 
$$ 2^{x_2} + b^{y_2} = 2 + (2^r h - 1)^{y_2} = 2^r h_1 + 1 = c^{z_2} = 2^{s + v_2(z_2)} k_1 + 1, 2 \nmid h_1, 2 \nmid k_1, 3 \mid h_1  \eqno{(4.18)}$$ 
where $v_2(n)$ is the greatest integer $t$ such that $2^t \mid n$.  From (4.18) we see that $r>s$, so that in (4.2) we must have $x_1 = s$.  

Now we apply Lemma 2.13.  Clearly $2< \frac{c^{z_2/2}}{4}$ (recall Lemma 2.9), so that applying Lemma 2.13 to (4.2), we find that we must have 
$$2^{x_1} = 2^s \ge \frac{c^{1/2}}{4} > \frac{b^{y_1/2}}{4} \ge \frac{b}{4} \ge \frac{3 \cdot 2^r - 1 }{4} \ge \frac{3 \cdot 2^{s+1} -1}{4} = 3 \cdot 2^{s-1} - \frac{1}{4} $$
 so that
$$ 2 \ge 3 - \frac{1}{2^{s+1}}, $$
which is false for all positive $s$.  

Thus we have $b \not\equiv 2 \bmod 3$.  So $b \equiv 1 \bmod 3$ and (4.6) holds.  
\end{proof}

We now prove Corollary 1.6:  

\begin{proof}
$N(a,b,c)=2$ follows from Theorem 1.4.  
If $z_1=z_2$, then, recalling Lemma 2.7 and noting that $c \ne 35$, $133$, or $259$, we must have  $b^{y_2} - 2^{x_1} = b^{y_1} - 2^{x_2} > 0$, so that Lemma 2.1 gives $2 \nmid x_1 - x_2$.  Consideration modulo 3 shows that this requires $2 \nmid y_1 - y_2$ and $b \equiv 2 \bmod 3$, contradicting Theorem 1.5.  So $z_1 \ne z_2$, so that by Theorem 1.4 $z_2>1$.        
\end{proof}

\section{The Unlikelihood of Counterexamples to Conjecture 1.7}  %5

We outline the algorithm used to justify Theorem 1.8.  

For a given prime value of $b$ and some small prime $p$ (or small prime power), we will consider all solutions $(x_1, y_1, x_2, y_2, z_2)$ to $2^{x_1} + b^{y_1} \equiv c \bmod p$ and $2^{x_2} + b^{y_2} \equiv c^{z_2} \bmod p$.  Note that the exponents are defined modulo $p-1$.  

When $b \equiv 13 \bmod 24$ and $c \equiv 5 \bmod 24$, we must have $x_1=2$.  From above $2 \mid y_1$, $2 \mid x_2$, $ 2 \nmid y_2$, $2 \nmid z_2>1$, and $z_2$ divides the class number of $\ratQ(-b)$. 

For a given $b \equiv 13 \bmod 24$ we find all $z_2>1$ with $z_2$ odd and dividing the class number.   Fix $b$ and $z_2$.  For a given prime $p$, we consider each $y_1 \bmod p-1$ with $y_1$ even.  Define $c \equiv 2^{2} + b^{y_1} \bmod p$; for each value of $x_2$ even and $y_2$ odd modulo $p-1$, we see if $ 2^{x_2} + b^{y_2} \equiv c^{z_2} \bmod p$.  If there is a solution, we add $(y_1,x_2,y_2)$ to a list of all possible solutions modulo this prime $p$.  We now consider another small prime (or prime power) $q$.  For each possible solution modulo $p$ we now see if there is a solution modulo $q$.    If we are fortunate, there are no solutions modulo $q$ that are consistent with a solution modulo $p$, in which case this choice of $b$ and $z_2$ cannot have any solutions. Otherwise, we add another prime $r$ and see if any solutions are consistent modulo $r$.  Often a given $b$, $z_2$ has no consistent solutions after checking a few primes (or prime powers).  In rare instances, the program required up to fourteen primes or prime powers to eliminate all possible solutions for a given $b$ and $z_2$, but we never needed primes exceeding 241.   

For $b \equiv 13 \bmod 24$ and $c \equiv 17 \bmod 24$, we must have $x_2=2$. The procedure is similar except that we now consider tuples $(x_1, y_1, y_2)$.  For $b \equiv 1 \bmod 24$, we cannot specify $x_1$ and $x_2$ so there are many more cases to check, but the same essential algorithm can be used. 

We used Maple{\textsuperscript{\textregistered}} for some preprocessing, then used Sage{\textsuperscript{\textregistered}} (in which we could access the 
Pari{\textsuperscript{\textregistered}} class number command) for the calculations. Total calculation time was about 100 hours.  See the third author's website for programs and details.  

In summary, we showed that for primes $b<10^9$, there are no solutions $(x_1, y_1, z_1, x_2, y_2, z_2)$ outside those listed in Conjecture 1.7.  Since $c>b^2$, we have $c>10^{18}$.  This concludes the demonstration of Theorem 1.8.  

In another direction, it is interesting to note that the case $(b,c) \equiv (13,17) \bmod 24$ requires the equation 
$$4+b^{y_2} = c^{z_2}, 2 \nmid y_2 > 1, 2 \nmid z_2>1. \eqno{(5.1)} $$
Note that if $y_2=1$, then $c=2^{x_1} + b^{y_1} > 4+b = c^{z_2}$, which is impossible; note also that $z_2=1$ is impossible by Corollary 1.6.  

The conditions in (5.1) are extremely unlikely even without the extra consideration of an additional solution $(x_1, y_1, z_1)$.


\begin{thebibliography}{1}


\bibitem{BB}
M. Bauer and M. A. Bennett, Applications of the hypergeometric method to the generalized Ramanujan-Nagell equation, {\it Ramanujan J.}, {\bf 6} (2002), 209--270. 

\bibitem{Be}
M. A. Bennett, On some exponential equations of S. S. Pillai, {\it Canadian Journal of Mathematics}, {\bf 53} no. 5 (2001), 897--922.

\bibitem{Be2} 
M. A. Bennett, Differences between perfect powers, {\it Canad. Math. Bull.}, {\bf 51} (2008), 337--347.


\bibitem{BS}   %3
M. A. Bennett and C. M. Skinner,
Ternary Diophantine Equations via Galois Representations and Modular Forms,
{\it Canad. J. Math.}, {\bf  56} (1) (2004), 23--54.




\bibitem{BSch}  %4
F. Beukers and H.P. Schlickewei,
The equation $x + y = 1$ in finitely generated groups,
{\it Acta Arith.},
{\bf 78} (1996),
189---199.


\bibitem{Cao}  %5
Z.-F. Cao,
\newblock On the diophantine equation $a^x + b^y = c^z$ {I},
\newblock {\em Chinese Sci. Bull.}, {\bf 31} (22) (1986), 1688--1690.
\newblock in Chinese.

\bibitem{Co}    %6  
J. H. E. Cohn, The Diophantine equation $x\sp 2+2\sp k=y\sp n$, {\it Arch. Math. (Basel)}, {\bf 59} no. 4 (1992), 341--344.


\bibitem{G}
A.~O. {Gel'fond},
\newblock Sur la divisibilit\'e de la diff\'erence des puissance de deux
  nombres entiers par une puissance d'un id\'eal premier,
\newblock {\em Mat. Sb.}, {\bf 7} (49) (1940), 7--25.

\bibitem{HL1}
Y.-Z. Hu and M.-H. Le,
\newblock An upper bound for the number of solutions of ternary purely
  exponential diophantine equations,
\newblock {\em J. Number Theory}, {\bf 183} (2018), 62--73.


\bibitem{HL2} 
Y.-Z. Hu and M.-H. Le, 
An upper bound for the number of solutions of ternary purely exponential diophantine equations II, 
{\it Publ. Math. Debrecen} {\bf 95} (2019),
335--354.



\bibitem{Le1}  %10
M.-H. Le,
\newblock On the diophantine equation $a^x + b^y = c^z$,
\newblock {\em J. Changchun Teachers College, Nat. Sci.}, {\bf 2} (1) (1985), 50--62.
\newblock in Chinese.


\bibitem{Le2} 
M.-H. Le, Some exponential Diophantine equations I: The equation $D_1 x^2 - D_2 y^2 = \lambda k^z$, {\em J. Number Theory},  {\bf 55}(2) (1955), 209--221.


\bibitem{Le}     %12
M.-H. Le,
\newblock On {Cohn}'s conjecture concerning the diophantine equation $x^2 + 2^m = y^n$,
\newblock {\em Archiv der Mathematik}, {\bf 78} (2002), 26--35.


\bibitem{LeScSt}
M. H. Le, R. Scott, R. Styer,  A survey on the ternary purely exponential Diophantine equation $a^x+b^y=c^z$, {\em Surv. Math. Appl.} (2019), 109--140.  


\bibitem{M}
K.~Mahler,
\newblock {Zur Approximation algebraischer Zahlen I: \"Uber den gr\"ossten
  Primteiler bin\"arer Formen},
\newblock {\em Math. Ann.}, {\bf 107} (1933), 691--730.

\bibitem{MP}
T. Miyazaki and I. Pink, Number of solutions to a special type of unit equation in two variables, arXiv:2006.15952.

\bibitem{MP2}
T. Miyazaki and I. Pink, Number of solutions to a special type of unit equation in two variables II, arXiv:2205.11217.


\bibitem{P}
S.S. Pillai, On the equation $2^x - 3^y = 2^X + 3^Y$, {\em Bull. Calcutta Math. Soc.}, {\bf 37} (1945), 18--20.


\bibitem{Sc}     %18
R. Scott, 
\newblock On the Equations $p^x-b^y = c$ and $a^x+b^y=c^z$, 
\newblock {\em Journal of Number Theory}, {\bf 44} no. 2 (1993), 153--165. 

\bibitem{ScSt1}    %19
R.~Scott and R.~Styer,
\newblock On $p^x - q^y = c$ and related three term exponential diophantine
  equations with prime bases,
\newblock {\em J. Number Theory}, {\bf 105} (2) (2004), 212--234.

\bibitem{ScSt7}     %20
R. Scott, R. Styer, Bennett's Pillai theorem with fractional bases and negative exponents allowed, {\em Journal de th\'eorie des nombres de Bordeaux}, {\bf 27}(1) (2015), 289--307.  


\bibitem{ScSt6}     %21
R. Scott, R. Styer, 
Number of solutions to $a^x + b^y = c^z$, 
{\em Publ. Math. Debrecen}
{\bf 88}
(2016),
 131--138.


\bibitem{StTi}
R.J. Stroeker, R. Tijdeman, Diophantine Equations, {\em Computational Methods in Number Theory, M.C. Tract 155}, Centre for Mathematics and Computer Science, Amsterdam (1982), 321--369.

\end{thebibliography}
\end{document}